\input amstex
\message{CatcodeOfAt=\the\catcode`\@}
\message{Abovedisplay=\the\abovedisplayskip}
\catcode`\@=11
\def\logo@{}
\catcode`\@=13
\tolerance=800
\NoBlackBoxes
\documentstyle{amsppt}
\NoRunningHeads

\def\Im{\operatorname{Im}}
\def\ctg{\operatorname{ctg}}
\def\tg{\operatorname{tg}}
\def\ch{\operatorname{ch}}
\def\sh{\operatorname{sh}}
\def\={\operatorname{=}}
\chardef\No=194

\TagsOnRight

\topmatter
\title
The asymptotical behavior of spectral function of one family of differential 
operators.
\endtitle
\author
Alexander S. Pechentsov and Anton Yu. Popov\\
Moscow State University, Russia
\endauthor

\endtopmatter

\document

In $L_2[0,+\infty )$ consider operator $l$, defined by differential
expression:
$$
ly(x)=-y^{\prime \prime }(x)+q(x)y(x)
$$

and boundary condition
$$
y(0)\cos \alpha +y^{\prime }(0)\sin \alpha =0.
$$

Assume that $\alpha \in R$, $q$ is a function continuous
on $[0,+\infty )$,
and having real values.
Denote by $\varphi (x,\lambda )$ and $\theta (x,\lambda )$
the solutions to equation $ly=\lambda y$ with initial conditions%
$$
\varphi (0,\lambda )=\sin \alpha ,\ \varphi ^{\prime }(0,\lambda )=-\cos
\alpha ,
\theta (0,\lambda )=\cos \alpha ,\ \theta ^{\prime }(0,\lambda )=\sin \alpha
.
$$

Functions $\varphi $,$\theta $ and $\rho $,$m$,$\hat f$ (we will introduce
them later) depend, naturally, on $\alpha $, but we will omit the argument $%
\alpha $, for not to overload the text with notations.

It is well known the Weyl theorem on representation of arbitrary function $%
f\in L_2[0,+\infty )$ as an integral along the spectrum of operator $L$.

\proclaim{Theorem [1-3]}  There exists such a non decreasing and bounded from
below on $\Bbb R$ function $\rho $ , for which  the following
statements take place:

\item{1.}
$
\exists \ \hat f(\lambda )\= \limits _{L_2(R,d\rho )}\lim\limits_{n%
\rightarrow \infty }\int\limits_0^nf(x)\varphi (x,\lambda
)dx=\int\limits_0^{+\infty }f(x)\varphi (x,\lambda )dx,
$
\newline
and the function $f$ is considered as a limit in $L_2[0,+\infty )$ of
''Fourier transform by measure $d\rho $ '' of $\hat f$ :

$$
f(x)\= \limits _{L_2(R^{+})}\lim \limits_{n\rightarrow \infty
}\int\limits_{-n}^n\hat f(\lambda )\varphi (x,\lambda )d\rho (\lambda
)=\int\limits_{-\infty }^{+\infty }\hat f(\lambda )\varphi (x,\lambda )d\rho
(\lambda ).
$$

\item{2.}  If $\sin \alpha \neq 0$, then for function $m(z)$ defined in
half-plane $\Im z>0$ by the equality
$$
m(z)=-\ctg \alpha +\int\limits_{-\infty }^{+\infty }\frac{%
d\rho (\lambda )}{z-\lambda },\tag1
$$
we have:
$$
\theta (x,z)+m(z)\varphi (x,z)\in L_2[0,+\infty )
\qquad
\forall z:\Im z>0.
$$

\item{3.}  For an arbitrary Borel set $E\subset R$ , not including the
spectrum of operator $L$, we have:$\int\limits_Ed\rho (\lambda )=0$.
\endproclaim

The statement 3 of  Theorem shows, that the measure $d\rho $ is
concentrated on the spectrum of operator $L$. The function $\rho (\lambda )$
is called the spectral function, and $m(z)$ is called the function of
Weyl-Titchmarsh of operator $L$. It is easy to see from the representation
(1) that the function $m(z)$ is analytical in the upper half-plane and
that if on certain interval $(a,b)\subset R$ there exists
$$
\lim \limits_{y\rightarrow +0}m(\lambda +iy)=m(\lambda )\in C(a,b),
$$
then $\rho (\lambda )$ has on this interval continuous derivative, which is
connected to $m(\lambda )$ by the correspondence:
$$
\rho ^{\prime }(\lambda )=\frac{-\Im m(\lambda )}{\pi} ,
\lambda \in (a,b). \tag2
$$
The investigation of properties of function $p(\lambda )$, and its
definition via the potential $q$ and the number $\alpha $, defining the
operator $L$, is, as a rule, a very difficult problem. There are well known the 
general
theorems of B.M.Levitan and V.A.Marchenko (see [4]---[6]
on evaluation of $\rho (\lambda )$ as $\lambda \rightarrow
\pm \infty $:
$$
\rho (\lambda )-\rho (-\infty )=o(\exp (-a\sqrt{|\lambda
|})),\ (\lambda \rightarrow -\infty )\ \forall a>0, \tag3
$$
$$
\rho (\lambda )=\frac{2\sqrt{\lambda }}{\pi \sin {}^2\alpha }+\rho (-\infty
)+\frac{\cos \alpha }{\sin {}^2\alpha }+o(1)\ (\lambda
\rightarrow +\infty ).
$$
(it is assumed, that the norming $\rho (0)=0$ is performed.)

The simplest case of $q(x)\equiv 0$ was studied by Titchmarsh in [7].
It turned out, that for $\lambda >0$
$$
\rho ^{\prime }(\lambda )=\frac{\sqrt{\lambda }}{\pi (\lambda \sin
{}^2\alpha +\cos {}^2\alpha )},
$$
and for $\lambda <0$ both of the cases  $\ctg\alpha >0$ and $\ctg
\alpha <0$ are strongly different:
$$
\rho _\alpha ^{\prime }(\lambda )\equiv 0,\ \lambda <0,\ \ctg \alpha
<0,
$$
$$
\rho _\alpha ^{\prime }(\lambda )=\frac{2\ctg \alpha }{\sin
{}^2\alpha }\delta (\lambda -\lambda _0),\ \lambda <0,\ \ctg \alpha
>0,\tag4
$$
where $\lambda _0=-\ctg^2\alpha $,$\ \delta $ is delta-function of
Dirac.

In this work we study the case of $\ctg\alpha >0$ and analyse the
behavior of the derivative of spectral function on $(-\infty ,0)$ for
potentials of special type.

The following problem of common form is interesting for our study:

Let $Q\in C(0,+\infty )$. Denominate as $\rho _\alpha (\varepsilon ,Q,\lambda 
)$
the spectral function of the operator, defined by differential expression:
$$
-y^{\prime \prime }(x)+\varepsilon Q(x)y(x)
$$
and boundary condition
$$
y(0)\cos \alpha +y^{\prime }(0)\sin \alpha =0,\ \ctg\alpha >0.
$$

Is it true that on $(-\infty,0 )$ the family of distributions $\rho
_\alpha ^{\prime }(\varepsilon ,Q,\lambda )$ converges as $\varepsilon
\rightarrow 0$ (in any ''reasonable'' sence) to the derivative of spectral
function (4) of non-perturbed operator $-y^{\prime \prime }$ ? The
positive solution of this problem would allow us to write for functions
$f \in L_2[0,+\infty )$ the ''approximate spectral expansion'' of the
following kind:
$$
f(x)=\int\limits_{\lambda _0-\beta }^{\lambda _0+\beta }+\int\limits_{-\beta
}^{+\infty }\hat f(\lambda )\varphi _\alpha (x,\lambda )d\rho _\alpha
(\varepsilon ,Q,\lambda )
$$
$(\beta =\beta (\varepsilon )\downarrow 0)$, omitting the integral on the
complement to small neighborhood of the spectrum of non-perturbed operator.
As we know, such investigations have not been performed before.

Titchmarsh in [7] have expressed through Gankel functions of the first
type the Weyl-Titchmarsh function $m(\lambda ,\varepsilon )$ of operator
$$
-y^{\prime \prime }-\varepsilon xy,\ \varepsilon >0,\ y(0)\cos \alpha
+y^{\prime }(0)\sin \alpha =0,\ \ctg\alpha >0.\tag5
$$

He proved, that for all $\lambda \in \Bbb R$

$$
m(\lambda ,\varepsilon )=\frac{H_{1/3}^{(1)}(A)\sin \alpha -\sqrt{%
\lambda }H_{-2/3}^{(1)}(A)\cos \alpha }{\text{$H_{1/3}^{(1)}(A)\cos \alpha +
\sqrt{\lambda }H_{-2/3}^{(1)}(A)\sin \alpha $}},\tag6
$$

where $A=\frac{2\lambda ^{3/2}}{3\varepsilon }$,$\lambda ^{3/2}=-i|\lambda
|^{3/2}$ with $\lambda <0$,$H_p^{(1)}$are Hankel functions of the first type.

Using the representation (6)  it was proved in \cite{6}  , that in certain
vicinity of the point $\lambda _{0\text{ }}$the functions $m(\lambda
,\varepsilon )$ have only one pole $\lambda _0(\varepsilon )$ with the 
asymptotics
$$
{\lambda _0 (\varepsilon )=\lambda _0-\frac{\varepsilon }{2}\tg
\alpha +O(\varepsilon ^2)}, \quad (\varepsilon \rightarrow +0),
$$
and its imaginary part satisfies the inequality:
$$
-\exp(-\frac{\ctg^3\alpha }{\varepsilon }) \leqslant \Im \lambda _0(
\varepsilon ) <0.
$$
But the behavior of $\rho ^{\prime}(\lambda ,\varepsilon )$ as
$\lambda \in (-\infty ,0)$ and $\varepsilon \rightarrow +0$ was
not studied by Titchmarsh. We could solve this problem in certain
sence. In spite of the existence of explicit formula for
$\rho ^{\prime} (\lambda, \varepsilon)$ (taking into account (2) we
have $\rho ^{\prime}(\lambda ,\varepsilon)=-\frac{\Im m(\lambda, 
\varepsilon)}{\pi} $,
due to the fact that from (6) it follows, that the function $m(\lambda, 
\varepsilon)$
is analytical in closed upper half-plane except one point $0$) this was
not easyly to realize. The matter is that in putting in (6) known
asymptotic series for $H_{1/3}^{(1)}(A)$ and $H_{-2/3}^{(1)}(A)$
as $A \rightarrow -i\infty $ (along negative part of imaginary axis),
then for $m(\lambda, \varepsilon)$ we obtain the asymptotic series, all
members of which are real. Therefore, mentioned series does not give 
information
concerning $\rho ^{\prime} (\lambda, \varepsilon)$.

We have obtained the following representation for function
$\rho ^{\prime} (\lambda, \varepsilon)$, which takes place
for any $\lambda < 0$ and $\varepsilon > 0$:
$$
\rho ^{\prime} (\lambda, \varepsilon)=
\frac
{\pi^{-1}\tau(\beta_{1,1}(a)\beta_{2,2}(a)+
\beta_{1,2}(a)\beta_{2,1}(a))}
{(\beta_{1,1}(a)\cos \alpha - \beta_{1,2}(a)\tau \sin \alpha)^2+
(\beta_{2,1}(a)\cos \alpha + \beta_{2,2}(a) \tau \sin \alpha )^2},
\tag7
$$
where $\tau = \sqrt{-\lambda}$, $a = |A| = \frac{2 \tau^2}{3 \varepsilon}$,
$$
\beta_{1,1}(a) = 1+ \frac{\sqrt{3}}{2} \frac{B^{-}_{1/3}(a)}
{\Omega _{1/3}(a)}, \quad
\beta_{1,2}(a)=\frac
{\Omega_{2/3}(a)}
{\Omega_{1/3}(a)}+
\frac{\sqrt{3}}{2}
\frac
{B^{-}_{2/3}(a)}
{\Omega_{1/3}(a)},
$$
$$
\beta_{2,1}(a)=
\frac
{B^{+}_{1/3}(a)}
{2\Omega_{1/3}(a)},
\quad
\beta_{2,2}(a)=
\frac
{B^{+}_{2/3}(a)}
{2\Omega_{1/3}(a)},
$$
$$
\Omega_{p}(a)=
\frac{1}{\pi}
\int\limits^{\pi}_{0}
\exp{(a \cos{t})\cos(pt)}dt,
$$
$$
B^{-}_{p}(a)=
\frac{1}{\pi}
\int\limits^{+\infty}_{0}
\exp{(-a \ch{t})} \sh(pt) dt,
$$
$$
B^{+}_{p}(a)=
\frac{1}{\pi}
\int\limits^{+\infty}_{0}
\exp(-a \ch{t}) \ch(pt) dt.
$$

The representation (7) allows us to prove the following
theorems 1-3:
\proclaim{Theorem 1}
For any $\varepsilon >0$  the following asymptotics takes place:
$$
\rho^{\prime}(\lambda, \varepsilon)
\sim
\frac
{\exp(-\frac{4}{3} \frac{|\lambda|^{3/2}}{\varepsilon})}
{\pi \sqrt{-\lambda} \sin^{2} \alpha},
\quad
(\lambda \rightarrow -\infty).
$$
\endproclaim
\proclaim{Theorem 2}
There exist positive constants $c_1$ and $c_2$,
effectively dependent on $\alpha$,
such that with $\varepsilon \in (0,c_1)$
take place evaluations:
$$
\rho^{\prime}(\lambda, \varepsilon)=
O(|\lambda|^{-1/2} \exp(-\frac{4|\lambda|^{3/2}}{3 \varepsilon})),
\quad
\lambda \in (-\infty, -2 \ctg^2 \alpha),
$$
$$
\rho^{\prime}(\lambda, \varepsilon)=
O(|\lambda|^{1/2} \exp(-\frac{4|\lambda|^{3/2}}{3 \varepsilon})),
\quad
\lambda \in (-\frac{1}{2} \ctg ^2 \alpha, -c_{2} \varepsilon ^{2/3}),
$$
$$
\rho^{\prime}(\lambda, \varepsilon)=
O(\varepsilon ^{1/3}),
\quad
\lambda \in [-c_2 \varepsilon^{2/3} ,0).
$$
Constants in symbols $O$ depend only on $\alpha $ and
are effective.
\endproclaim

Denote $T(\lambda, \varepsilon)= \beta_{1,1}(a) \cos \alpha
 - \beta_{1,2}(a) \tau \sin \alpha $. The answer to the
question on the behavior of $\rho ^{\prime}(\lambda, \varepsilon)$
on the "critical" segment
$ I=[-2 \ctg ^2 \alpha , -\frac{1}{2} \ctg ^2 \alpha ]$
gives the following

\proclaim{Theorem 3}
There exists such a positive constant $c_3$, depending effectively on $\alpha $
, that for any $\varepsilon \in (0,c_3)$  the following statements take place:

1) Function $T(\lambda, \varepsilon)$ has on segment $I$ unique
zero $\lambda_1 (\varepsilon)$ with asymptotics:
$$
\lambda_1 (\varepsilon)=
\lambda_0 - \frac{\varepsilon}{2} \tg \alpha + O(\varepsilon ^2).
$$

2) For any $d>0$,$\lambda \in I\backslash I_{d,\varepsilon}
\qquad
(I_{d,\varepsilon}=[\lambda_1(\varepsilon)-d, \lambda_1(\varepsilon)+d])$
 the evaluation takes place:
$$
\rho^{\prime}(\lambda, \varepsilon )=
d^{-2} \exp(-\frac{4|\lambda|^{3/2}}{3 \varepsilon }).
$$

3) If $\exp (-\frac{\ctg ^3 \alpha }{3 \varepsilon }) \le d \le \varepsilon^2$,
then
$$
\int\limits_{I_{d, \varepsilon }} d \rho (\lambda) =
\frac
{2 \ctg \alpha}
{\sin ^2 \alpha }+
O(\varepsilon)
$$
\endproclaim
\proclaim{Corollary}
In space $C^{*}(I)$ the family of functions
$\rho ^{\prime} (\lambda, \varepsilon)$
converges weakly to
$\frac{2 \ctg \alpha }{\sin ^2 \alpha } \delta (\lambda -
\lambda_0 )$ as $\varepsilon \rightarrow +0$.
\endproclaim

\bigskip
{\smc       Alexander S. Pechentsov and Anton Yu. Popov,
Moscow State University, Russia.}

\enddocument